\documentclass[12pt]{article}
\usepackage{latexsym,amsmath,amssymb}

\newif\ifdviwin

\usepackage[latin1]{inputenc}
\usepackage[english]{babel}
\usepackage{indentfirst}
\usepackage[mathscr]{eucal}
\usepackage{amssymb,amsmath,amsfonts}
\usepackage{fancybox,fancyhdr}
\usepackage{graphicx}

\newif\ifdviwin

\dviwintrue

\def\cS{\mathcal{S}}

\def\cD{\mathcal{D}}

\def\cM{\mathcal{M}}
\def\cG{\mathcal{G}}
\def\cF{\mathcal{F}}

\let\8=\infty \let\0=\emptyset

\let\landa=\lambda
\let\alfa=\alpha

\let\parc=\partial

\def\vart{\vartheta}
\def\landa{\lambda}
\def\Landa{\Lambda}

\def\flecha{\rightarrow}

\def\cte.{\mathop{\rm cte.}\nolimits}

\def\Re{\mathop{\rm Re }\nolimits} \def\Im{\mathop{\rm Im }\nolimits}

\def\N{\mathbb{N}}

\def\R{\mathbb{R}}

\def\C{\mathbb{C}}
\def\D{\mathbb{D}}

\def\S{\mathbb{S}}

 \newtheorem{defi}{Definition}
 \newtheorem{teo}[defi]{Theorem}
 
 \newtheorem{cor}[defi]{Corollary}
 \newtheorem{lem}[defi]{Lemma}
 
 \newtheorem{remark}[defi]{Remark}

 \newenvironment{proof}{\rm \trivlist \item[\hskip \labelsep{\it
      Proof}:]}{\par\nopagebreak \hfill $\Box$ \endtrivlist}

 \newenvironment{proof1}{\rm \trivlist \item[\hskip \labelsep{\it
      Proof of Theorem \ref{main}}:]}{\par\nopagebreak \hfill $\Box$ \endtrivlist}

\numberwithin{equation}{section} 

\begin{document}
\mbox{}\vspace{0.4cm}\mbox{} 

\begin{center}
\rule{15.2cm}{1.5pt}\vspace{0.5cm} 

{\Large \bf The space of solutions to the Hessian one equation 
\\[0.3cm] in the finitely punctured plane}\\ \vspace{0.5cm} {\large José A.
Gálvez$\mbox{}^a$, Antonio Martínez$\mbox{}^b$ and Pablo 
Mira$\mbox{}^c$}\\ \vspace{0.3cm} \rule{15.2cm}{1.5pt} 
\end{center}
  \vspace{1cm}
$\mbox{}^a$, $\mbox{}^b$ Departamento de Geometría y Topología, Universidad de Granada,
E-18071 Granada, Spain. \\ e-mail: jagalvez@ugr.es ; amartine@ugr.es 
\vspace{0.2cm} 

\noindent $\mbox{}^c$ Departamento de Matemática Aplicada y Estadística, 
Universidad Politécnica de Cartagena, E-30203 Cartagena, Murcia, Spain. \\ 
e-mail: pablo.mira@upct.es \vspace{0.2cm} 

\noindent Date: September 10, 2004  

\vspace{0.3cm} 

 \begin{abstract}
We construct the space of solutions to the elliptic Monge-Ampère equation 
${\rm det} (D^2 \phi)=1$ in the plane $\R^2$ with $n$ points removed.
We show that, modulo equiaffine transformations and for $n>1$, this space 
can be seen as an open subset of $\R^{3n-4}$, where the coordinates are 
described by the conformal equivalence classes of once punctured bounded 
domains in $\C$ of connectivity $n-1$. This approach actually provides a 
constructive procedure that recovers all such solutions to the 
Monge-Ampère equation, and generalizes a theorem by K. Jörgens. 
 \end{abstract}



\section{Introduction}
A celebrated result by K. Jörgens \cite{Jor1} states that all solutions to 
the elliptic Monge-Ampère equation 
\begin{equation}\label{monge}
\phi_{xx} \phi_{yy} - \phi_{xy}^2 =1 
\end{equation}
which are globally defined on $\R^2$ are quadratic polynomials. This 
theorem has motivated a large amount of works, which essentially follow 
two research lines. One of them deals with the extension of the result to 
more general classes of Monge-Ampère equations (see for instance 
\cite{BCGJ,Cal,CaLi2,ChYa,ChWa,Pog,TrWa1,TrWa2}). The other one was 
initiated by K. Jörgens in 1955, and concerns the validity of the theorem 
in large proper domains of $\R^2$, and the study of the new solutions that 
may arise. Some results of this type can be found in 
\cite{Jor2,CaLi1,FMM1,FMM2,TrWa1}. In particular, in \cite{Jor2} all the 
solutions to \eqref{monge} globally defined in $\R^2 \setminus \{(0,0)\} $ 
were obtained. 

In the present work we follow this last direction, and plan to determine 
all the solutions to \eqref{monge} in $\R^2$ with $n$ points removed. In 
this way, our result generalizes the two mentioned theorems by Jörgens, 
and suggests a new line of inquiry into the theory of Monge-Ampère 
equations. 

For stating such a general description, some remarks have to be made. 
First, note that if $\phi$ is a solution to \eqref{monge} in a finitely 
punctured plane, the punctures correspond then to isolated singularities 
of $\phi$. We shall assume without loss of generality that isolated 
singularities are \emph{non-removable}, i.e. $\phi$ cannot be 
$\mathcal{C}^1$-extended across the singularity. A well-known criterion to determine
whether an isolated singularity of a solution to \eqref{monge} is 
removable is found in \cite{Jor2}. We shall provide a different one in 
Section 3, in terms of the underlying conformal structure of the solution 
$\phi$.

Another remark concerns non-unicity of the solutions to \eqref{monge}. Let 
$\mathcal{E} $ be the group of \emph{equiaffine transformations} of $\R^3$
whose differentials fix the $x_3$-direction, i.e. the set of maps $\Phi 
:\R^3 \flecha \R^3$ of the form 
 \begin{equation}\label{equi}
 \Phi \left(\def\arraystretch{1.2}\begin{array}{c} x_{1} \\
x_{2}\\ x_{3} 
\end{array} \right) =
\left(\def\arraystretch{1.2}\begin{array}{ccc} a_{11} & a_{12} & 0 \\ 
a_{21} & a_{22} & 0 \\ a_{31} & a_{32} & 1 
\end{array} \right) \left(\def\arraystretch{1.2}\begin{array}{c} x_{1} \\
x_{2}\\ x_{3} 
\end{array} \right) + \left(\def\arraystretch{1.2}\begin{array}{c} b_{1} \\
b_{2}\\ b_{3} 
\end{array} \right), \hspace{0.5cm} a_{11} a_{22} - a_{21} a_{12} =1.
 \end{equation}
Then, if $\phi(x,y)$ solves \eqref{monge} and $\Phi \in \mathcal{E}$, the 
function $\phi' (x',y')$ given by $\Phi (x,y,\phi (x,y))= 
(x',y',\phi'(x',y'))$ is a new solution to \eqref{monge} in terms of the 
variables $x'$, $y'$. Therefore, we shall be interested in obtaining the 
classification result for solutions to \eqref{monge} \emph{modulo 
equiaffine transformations}. In other words, two solutions to 
\eqref{monge} differing only by an element of $\mathcal{E}$ will be 
considered to be equivalent. 

Keeping these facts in mind, the classification by Jörgens of the 
solutions to \eqref{monge} in the once punctured plane is formulated as 
follows: \emph{the unique solutions to \eqref{monge} in 
$\R^2$ and in $\R^2 \setminus \{(0,0)\}$ are respectively, up to equiaffine transformations, the
quadratic polynomial $\phi(x,y)= \left(x^2 +y^2 \right)/2$ and the 
rotational example} 
$$\phi(x,y)= \frac{1}{2}\left(\sqrt{x^2 + y^2} \sqrt{1 + x^2 + y^2} + \sinh^{-1} \left(\sqrt{x^2 + y^2} \right)
\right).$$ So, the following question arises naturally: given $n>1$, is 
there a unique (up to equiaffine transformations) solution to 
\eqref{monge} defined in an $n$-times punctured plane, and so that the 
punctures are non-removable isolated singularities? In Section 3 we shall 
analyze to what extent the answer to this question is \emph{yes}. 

The paper is organized as follows. In Section 2 we review the theory of 
solutions to \eqref{monge}, and its equivalence to the theory of locally 
convex embedded parabolic affine spheres in the affine $3$-space 
$\R^3$. In particular, we shall recall a conformal representation
available for these surfaces, and indicate some basic properties of 
parabolic affine spheres in terms of the holomorphic data. The main result 
in this Section is the construction of solutions to \eqref{monge} in a 
punctured plane $\R^2 \setminus \{p_1, \dots, p_n\}$ so that the $p_j$'s 
are non-removable isolated singularities of such solutions. The 
construction relies on the conformal representation of parabolic affine 
spheres, and the existence of two classical holomorphic functions with the 
fundamental property of mapping bijectively a once punctured bounded 
domain in $\C$ with $n$ boundary components onto 
$\C$ with $n$ vertical (resp. horizontal) slits.

Section 3 contains the main theorem of this work, which characterizes the 
examples constructed in Section 2 as the only solutions (up to equiaffine 
transformations) to \eqref{monge} that are globally defined over a 
finitely punctured plane. This theorem provides an explicit one-to-one 
correspondence between the conformal equivalence classes of $\C$ with $n$ 
disks removed, and the space of solutions to \eqref{monge} in a finitely 
punctured plane with exactly $n$ non-removable isolated singularities, 
modulo equiaffine transformations. Particularly, for $n>1$, this quotient 
space can be naturally identified with an open subset of $\R^{3n-4}$. 

Finally, in Section 4 we show that our main theorem for the particular 
choice $n=1$ recovers the Jörgens one, and we determine explicitly in 
terms of theta functions all the solutions to \eqref{monge} in the twice 
punctured plane. This indicates that all the solutions we have constructed 
are highly non-trivial. 

\section{Construction of the canonical examples}

Let $\phi(x,y):D\subseteq \R^2\flecha \R$ be a solution to \eqref{monge} 
on a planar domain $D$. We shall assume without loss of generality that 
$\phi_{xx}$ is always positive. Then the subset
$\mathcal{S}_{\phi} =\{(x,y,\phi(x,y)): (x,y)\in D\} \subset \R^3$
describes an embedded locally convex \emph{parabolic affine sphere} in the 
affine $3$-space $\R^3$ with affine normal vector $\xi =(0,0,1)$ and 
affine metric $ds^2 = \phi_{xx}\, dx^2 + 2 \phi_{xy}\, dx dy + \phi_{yy}\, 
dy^2$ (see for instance \cite{FMM1,FMM2} for more information). 
Conversely, any embedded locally convex parabolic affine sphere with 
affine normal vector $\xi=(0,0,1)$ is locally the graph over a domain in 
the $x,y$-plane of a solution to \eqref{monge}. Along this paper, by a 
parabolic affine sphere we will always mean a locally convex one. 

The affine metric $ds^2$ of $\mathcal{S}_{\phi}$ induces a Riemann surface 
structure on $\mathcal{S}_{\phi}$. We shall call this Riemann surface 
structure the \emph{underlying conformal structure} of 
$\phi(x,y)$.

The solutions to \eqref{monge} can be explicitly recovered in terms of 
holomorphic data with respect to this underlying conformal structure. We 
expose here the global version of this fact in \cite{FMM2}, for the 
general case of immersed parabolic affine spheres. 

\begin{teo}\label{reptem}
Let $\Sigma$ be a simply connected Riemann surface, and let 
$F,G:\Sigma\flecha \C$ be two holomorphic functions for which $dG\neq 0$ and $|dF|<|dG|$ hold
everywhere. Then the map $\psi:\Sigma\flecha \R^3$ given by 
\begin{equation}\label{rep}
\psi = \frac{1}{2} \left( G +\overline{F}, \frac{1}{4} \left\{ |G|^2 - 
|F|^2 + 2\Re (G F) \right\} - \Re \int F dG \right) 
\end{equation}
is a parabolic affine sphere with affine normal vector $\xi=(0,0,1)$ and 
affine metric $ds^2 =\frac{1}{4} \left( |dG|^2 - |dF|^2 \right)$. 

Conversely, any parabolic affine sphere can be recovered in this way with 
respect to the conformal structure induced by its affine metric. 
\end{teo}


In order to construct examples of solutions to \eqref{monge} in the 
$n$-times punctured plane, we need to recall first a
classical description of two special holomorphic functions, which we will 
denote by $p(z), q(z)$. 

Let $\Omega\subset \C$ be a bounded domain whose boundary is made up by 
$n\geq 1$ analytic Jordan curves $C_1, \dots, C_n$, where we will let $C_n$ be the exterior
component, and take $z_0 \in \Omega$. Let $\cG (z,z_0)$ be the Green 
function of 
$\Omega$ with respect to $z_0$, given by $\cG(z,z_0)= \mathcal{D} (z,z_0)
- \log |z-z_0|$, where $\cD$ is the solution to the Dirichlet problem in 
$\Omega$ with boundary conditions
$\cD (\zeta) = \log |\zeta - z_0| $, $\zeta \in C_k$. Then the map
 \begin{equation}\label{aru}
 u (z)=  \frac{\parc \cG (z,z_0)}{\parc z_0} + \frac{\parc \cG (z,z_0)}{\parc \overline{z_0}}
 \end{equation}
is harmonic in $\Omega \setminus \{z_0\}$, and vanishes on $\parc \Omega$. 
Now let $\omega_j$ denote the \emph{harmonic measure} of $C_j$ with 
respect to 
$\Omega$, characterized by the conditions
 \[ \left\{\def\arraystretch{1.2} \begin{array}{lll} \Delta \omega_j & = & 0 \hspace{0.5cm} \text{ on } \Omega, \\ \omega_j (C_j) &
 = & 1, \\ \omega_j (C_k) & = & 0 \hspace{0.5cm} \text{ if } k \neq j.
  \end{array} \right. \] We shall denote by $\alfa_{j,k}$ the period of the
conjugate harmonic function of $\omega_j$ along $C_k$: \[ \alfa_{j,k} = 
\int_{C_k} \frac{\parc \omega_j}{\parc {\bf n}} \ ds,\] where ${\bf n}$ is 
the unit normal to the curve $C_k$. Let finally $A_k \in \R$ denote the 
period of the conjugate harmonic function of $u(z)$ along $C_k$. Then, if 
$\landa_1, \dots, \landa_{n-1}$ are real numbers, the harmonic function
 \begin{equation}\label{truqui}
  u + \landa_1 \omega_1 +\dots + \landa_{n-1} \, \omega_{n-1} :
\Omega\setminus \{z_0 \} \flecha \R 
 \end{equation}
has a well defined conjugate harmonic function on $\Omega \setminus 
\{z_0\}$ if and only if the 
$\landa_i$'s solve the following linear system:

\[ \left\{\def\arraystretch{1.2} \begin{array}{ccc}
\landa_1 \, \alfa_{1,1} + \cdots + \landa_{n-1}\, \alfa_{n-1,1} & = & -A_1 
\\ \vdots & & \vdots \\ \landa_1 \, \alfa_{1,n-1} + \cdots + \landa_{n-1}
\, \alfa_{n-1,n-1} & = & -A_{n-1} 
  \end{array} \right. \]

This system is classically known to have a unique solution (see \cite{Ahl} 
for instance), and thus there is a holomorphic function $p(z): 
\Omega\setminus \{z_0\} \flecha \C$ whose real part is the above harmonic 
function, i.e. 
$\Re p(z)= u(z) + \sum_{i=1}^{n-1} \landa_i \, \omega_i (z)$. We must 
remark that $p(z)$ is well defined up to pure imaginary additive 
constants. 

In addition, by the construction of $u(z)$ it becomes clear that $p(z)$ 
has a pole of order one at $z_0$, with residue $1$. Moreover, $\Re p(z) = 
\landa_k \in \R$ along $C_k$ for every $k\in \{1, \dots, n\} $. 

Actually, $p(z)$ is characterized by these conditions: it is the unique 
(up to pure imaginary additive constants) holomorphic map in $\Omega 
\setminus \{z_0\}$ with a simple pole of residue $1$ at $z_0$ whose real 
part is constant along each boundary component in $\parc \Omega$. Finally, 
$p(z)$ is a conformal equivalence between $\Omega \setminus \{z_0\}$ and 
its image, which is a vertical slit domain in $\C$ \cite{Ahl}. 

The analogous process to the above one, but this time starting with 
\begin{equation}\label{ardo}
 v (z)= \frac{\parc \cG (z,z_0)}{\parc
 z_0} - \frac{\parc \cG (z,z_0)}{\parc \overline{z_0}},
 \end{equation}
produces a holomorphic function $q(z):\Omega \setminus \{z_0\} \flecha \C$ 
with a simple pole of residue $1$ at $z_0$, and such that its imaginary 
part $\Im q(z)$ is constant along each boundary curve $C_k$. Again, $q(z)$ 
is defined up to real additive constants, and it is the unique holomorphic 
map in $\Omega\setminus \{z_0\}$ satisfying these properties. Further, it 
maps $\Omega\setminus \{z_0\} $ bijectively onto a horizontal slit domain 
in $\C$. 

\begin{teo}\label{ejes}
For any $n\in \N$, there exists a regular solution $\phi(x,y)$ to 
\eqref{monge} globally defined in $\R^2$ minus $n$ points, and such that 
the punctures are non-removable isolated singularities of $\phi(x,y)$. 
\end{teo}
\begin{proof}
We start with a complex domain $\Omega$ as above, with a distinguished 
point $z_0\in \Omega$. Let then $p(z),q(z)$ be the holomorphic functions 
constructed following the previous procedure. Now, define 
\begin{equation}\label{fg}
G(z)= p(z) + q(z), \hspace{1cm} F(z)= p(z) -q(z). 
\end{equation}
It then follows readily that $(G + \overline{F})|_{C_k}$ is constant for 
all $k\in \{1,\dots, n\} $. Differentiating now this expression we infer 
that $|dF / dG| =1$ at every point in $\parc \Omega$. On the other hand, 
since 
$p(z),q(z)$ have simple poles of equal residues at $z_0$, we obtain that
$|dF /dG| (z_0)=0$. Hence, $|dF / dG |<1$ in all
$\Omega$. In addition, as $p,q$ are holomorphic bijections between $\Omega\setminus{\{z_0\}}$
and a certain slit domain in $\C$, we get that $dF,dG$ cannot vanish 
simultaneously. Particularly, $dG$ never vanishes by $|dF /dG|<1$. 

Therefore, by Theorem \ref{reptem} the holomorphic functions $F,G$ define 
a parabolic affine sphere $$\psi : \widetilde{\Omega \setminus \{z_0\} } 
\flecha \R^3,$$ where $\widetilde{\Omega \setminus \{z_0\}} $ is the 
conformal universal covering of $\Omega \setminus \{z_0\}$. The affine 
sphere $\psi$ is actually well defined on $\Omega \setminus \{z_0\}$. To 
see this we have to check that 
  \begin{equation}\label{perio}
\Re \int_{\Gamma} F dG =0 
  \end{equation}
along any loop $\Gamma$ in $\Omega \setminus \{z_0\}$. First, observe that 
as $p(z)$ and $q(z)$ map each $C_k$ into a piece of a straight line, the 
functions $F,G$ can be holomorphically extended across each $C_k$ by 
Schwarzian reflection. In this way, we may extend 
$\Omega \setminus \{z_0\}$ slightly so that the curves $C_k$ are interior to
the extended domain, and 
$F,G$ are holomorphic there. But in this larger domain, the curves $C_1, \dots, C_n$
provide a basis of its first homology group. So, we just need to check 
condition \eqref{perio} for $\Gamma = C_1, \dots, C_n$. To do so, we first 
observe that by \eqref{fg} it holds for all $k\in \{1,\dots, n\}$ 
$$\int_{C_k} F dG = \int_{C_k} \left( p\,  dq - q \, dp \right).$$ In addition, $$0 = \Re
\int_{C_k} d \left( p \bar{q} \right) = - \Re \int_{C_k} \left( p\, dq - q 
\, dp \right) +2 \int_{C_k} \left( \Im q \, \Im p' + \Re p \, \Re 
q'\right) dz.$$ As both $\Re p(z)$ and $\Im q(z)$ are constant along 
$C_k$, we obtain that $\Re \int_{C_k} F dG =0 $ for
all the $C_k$'s, as we wished. 

Finally, we show that $\psi=(\psi_1,\psi_2,\psi_3):\Omega \setminus 
\{z_0\} \flecha \R^3 
$ is actually an embedding. For this, observe first of all that $F$ is
holomorphic at $z_0$, and $G$ has a simple pole at $z_0$. Thus, 
$(\psi_1,\psi_2)$ is a homeomorphism from a small enough punctured disk about $z_0$
onto the exterior of a Jordan curve in $\R^2$. 

Next, observe that $\cF := (\psi_1, \psi_2):\Omega \setminus 
\{z_0\}\flecha \R^2$ is well defined and continuous on the once punctured 
topological sphere $E$ which results when each curve $C_k$ in 
$\overline{\Omega} \setminus \{z_0\}$ is identified with a single point.
This happens because by \eqref{rep} we have that $\mathcal{F} (C_j)=p_j 
\in \R^2$ for every $j\in \{1,\dots, n\} $ and for some points 
$p_1,\dots, p_n \in \R^2$.

By the behaviour of $\cF$ about $z_0$, there is a topological punctured 
disk $\mathcal{D} \subset E$ with puncture at $z_0$ such that 
$\cF $ maps $\cD$ bijectively onto the exterior of a closed disk $B$ in $\R^2$.
Thus, $\cF (E\setminus \cD)$ must fill the whole $\overline{B}$. As a 
consequence, $(\psi_1,\psi_2):\overline{\Omega}\setminus \{z_0\} \flecha 
\R^2$ is onto. 

Now, assume the existence of two points $z_1,z_2\in \Omega \setminus 
\{z_0\}$ with $\cF (z_1)=\cF(z_2) =a\in \R^2$, and let $C$ denote a 
divergent curve in $\R^2\setminus \{p_1,\dots, p_n\}$ with endpoint $a$ 
(that may belong to $\{p_1,\dots, p_n\}$). As $\cF$ is a local 
homeomorphism, there exist two paths $\gamma_1, \gamma_2$ in 
$\Omega$ with endpoints $z_1,z_0$ and $z_2,z_0$ respectively, such that
$\cF (\gamma_i)=C$. But $\cF$ is a homeomorphism around $z_0$, and this
implies that $\gamma_1, \gamma_2$ must coincide from a \emph{first point} 
to the endpoint $z_0$. This contradicts that $\cF$ is a local 
homeomorphism, unless $z_1=z_2$. Hence, 
$\psi :\Omega \setminus \{z_0\} \flecha \R^3$ is one-to-one.

Finally, we need to check that $\cF (\Omega\setminus \{z_0\}) \subseteq 
\R^2 \setminus \{p_1,\dots, p_n\} $. If $\cF (\widetilde{z})= p_j$ for 
some 
$\widetilde{z}\in \Omega$, we may take disjoints neighbourhoods $A_1,A_2$ of $\widetilde{z}$ and $C_j$
in $\overline{\Omega} \setminus \{z_0\}$, respectively. There exist then 
points $\zeta_i\in A_i$ with $\cF(\zeta_1)=\cF(\zeta_2) \in \R^2 \setminus 
\{p_1, \dots, p_n\} $. Thus both points lie in $\Omega$, which is 
impossible. 

Summing up, we have an embedded parabolic affine sphere $\mathcal{S}=\psi 
(\Omega \setminus \{z_0\}) \subset \R^3$ such that $(\psi_1, \psi_2 
):\Omega \setminus \{z_0\} \flecha \R^2 \setminus \{p_1, \dots, p_n\} $ is 
a global homeomorphism. This produces a regular solution to \eqref{monge} 
globally defined in $\R^2 \setminus \{p_1, \dots, p_n\} $ and with 
isolated singularities exactly at the points $p_1, \dots, p_n$. All these 
isolated singularities are non-removable by Lemma \ref{remo}. 

\end{proof}

\section{The classification theorem}

This Section is devoted to show that any solution to \eqref{monge} in a 
finitely punctured plane must be equiaffinely equivalent to one of the 
examples of Theorem \ref{ejes}. In order to do so, it is convenient to 
introduce the following terminology. By definition, a \emph{global 
solution} to \eqref{monge} in $\R^2 \setminus \{p_1,\dots, p_n\} $ will be 
a $\mathcal{C}^2$ function $\phi (x,y): \R^2 \setminus \{p_1, \dots, p_n\} 
\flecha \R$ verifying \eqref{monge} and that is not $\mathcal{C}^1$ at the 
$p_j$'s. We remark (see \cite{Jor2}) that if a solution to \eqref{monge}
extends as a $\mathcal{C}^1$ function across an isolated singularity, then 
it actually extends analytically across it. 

With this, we shall prove 

\begin{teo}\label{main}
Any global solution $\phi (x,y):\R^2 \setminus \{p_1, \dots, p_n\} \flecha 
\R$ to \eqref{monge} is one of the examples in Theorem \ref{ejes} up to 
equiaffine transformations. 

In particular, there exists an explicit bijective correspondence between 
the conformal equivalence classes of once punctured bounded domains in 
$\C$ of connectivity $n-1$, and the space of global solutions to \eqref{monge}
with $n$ punctures, modulo equiaffine transformations. 
\end{teo}

Let $\cM_n$ be the quotient space of global solutions to \eqref{monge} 
with $n>1$ punctures, modulo equiaffine transformations. As an interesting 
consequence of Theorem \ref{main} we can endow $\cM_n$ with a finite 
dimensional analytic manifold structure: 

\begin{cor}\label{space}
The space $\cM_n$ can be naturally identified with an open subset of 
$\R^{3n-4}$.
\end{cor}
\begin{proof}
Let $\Omega \subset \C$ be a bounded domain in $\C$ of connectivity 
$n-1$, and take $z_0 \in \Omega$. There exists then a Möbius
transformation of the Riemann sphere $\C \cup \{\8\} $ taking $\Omega 
\setminus \{z_0\}$ to a new complex domain $\Landa\subset \C$ consisting 
of 
$\C$ with the interior of $n$ Jordan curves removed. Thus, by Koebe's
uniformization theorem, $\Landa$ is conformally equivalent to $\C$ with 
$n$ disjoint disks removed. Moreover, these disks $D_1,\dots, D_n$ are
uniquely determined if we assume the following two restrictions: $D_1 =\D$ 
is the unit disk, and 
$D_2$ has its center in the positive real line.

Let now $\mathcal{B}_n$ denote the space of conformal equivalence classes 
of once punctured bounded domains in $\C$ of connectivity $n-1$. The above 
comments show that $\mathcal{B}_n $ can be canonically identified with the 
open subset in $(1,+\8)\times \C^{n-2} \times (\R^+)^{n-1} \subset 
\R^{3n-4}$ consisting of those points $a = (c_1,\dots, c_{n-1}, r_1, 
\dots, r_{n-1})$ for which the open disks of centers $c_j$ and associated 
radii $r_j$ (with $c_1 \in (1,+\8)$) together with the unit disk $\D$ are 
pairwise disjoint. Once here, the proof follows immediately from Theorem 
\ref{main}. 
\end{proof}

Before coming to the proof of Theorem \ref{main}, we shall give a 
criterion to determine when an isolated singularity of a solution to 
\eqref{monge} can be removed. Our perspective here is completely different 
from the usual results on removing isolated singularities of Monge-Ampère 
equations that can be found in the literature. However this description is 
to some extent implicit in the result by Jörgens \cite{Jor2} on removing 
isolated singularities of \eqref{monge}. 

\begin{lem}\label{remo}
An isolated singularity of a solution $\phi(x,y)$ to \eqref{monge} is 
removable if and only if the underlying conformal structure of $\phi(x,y)$ 
around the singularity is that of a punctured disk. 
\end{lem}
\begin{proof}
Let $\phi(x,y)$ be a solution to \eqref{monge} in a once punctured 
topological disk $\mathcal{U}\setminus \{(x_0,y_0)\}\subset \R^2$. Let 
$\zeta :\Omega \subset \C\flecha \mathcal{U}$ be a global holomorphic
parameter with respect to the underlying conformal structure of 
$\phi(x,y)$. Then, shrinking $\mathcal{U}$ if necessary,
$\Omega$ is biholomorphic to a punctured disk or an annulus

Suppose that $(x_0,y_0)$ is a removable isolated singularity of 
$\phi$, i.e. $\phi$ extends as a $\mathcal{C}^1$ function across
$(x_0,y_0)$. As the holomorphic data verify (see \cite{FMM1} for instance)
$G+F = 2x +2 i \phi_y$, $G-F= 2 \phi_x + 2 i y $, both $F,G$ have a well defined value at the
singularity. This implies that we must have the punctured disk conformal 
type (otherwise, $F,G$ would be constant).

Now, suppose that the underlying conformal structure is that of a 
punctured disk, say $\Omega =\D^*=\{z :0<|z|<1\}$, and let $F,G$ be the 
holomorphic data of the conformal representation. As was proved in 
\cite{FMM2}, it must hold on $\D^*$ that $dx^2 + dy^2 \leq |dG|^2$. Hence 
$G$ has at most a pole at $0$, and $dG(0) \neq 0$. In addition, as $|dF
/dG|<1$ on $\D^*$, we conclude that $F$ has at most a pole at $0$. But 
now, if $F$ or $G$ had a pole at $0$, the same would happen to $G+F$ or 
$G-F$. This is not possible, by \eqref{rep}. So $F,G$ extend
holomorphically to $0$. Finally 
$|dF /dG|(0) \neq 1$, because $|dF /dG|<1$ on $\D^*$.
Therefore, the conformal representation assures that the parabolic affine 
sphere extends regularly across $0$, i.e. $(x_0,y_0)$ is a removable 
singularity of $\phi(x,y)$. 
\end{proof}

\begin{remark}
Lemma \ref{remo} indicates that if $\phi(x,y)$ is a solution to 
\eqref{monge} in a punctured neighbourhood of a non-removable isolated 
singularity $(x_0,y_0)$, it has the underlying conformal structure of an 
annulus $A_r = \{z: 1<|z|<r\} $. Moreover, it is deduced from the 
conformal representation that $(G + \overline{F} )|_{\S^1} = x_0 +i y_0$, 
where $\S^1 =\{z: |z|=1\} $. So, both $G$ and $F$ can be analytically 
continued across $\S^1$ by Schwarzian reflection. In particular, 
$\phi(x,y)$ extends continuously to $(x_0,y_0)$, although it is not
$\mathcal{C}^1$ at this point.
\end{remark}

\begin{proof1}
Let $\phi(x,y):\R^2\setminus\{p_1,\dots, p_n\} \flecha \R$ be a global 
solution to \eqref{monge}, and let $\Sigma$ denote its underlying Riemann 
surface structure. By Lemma \ref{remo}, the underlying conformal type of 
$\phi(x,y)$ about any of the isolated singularities is that of an annulus.
Moreover, it was proved in \cite{FMM2} that the conformal type of 
$\phi(x,y)$ in the exterior of a sufficiently large disk $x^2 +y^2 >R^2$
is that of a punctured disk. Putting together these conditions and 
applying Koebe's uniformization theorem we deduce that $\Sigma$ is 
conformally equivalent to $\C$ with $n$ disjoint disks removed. We shall 
denote by $\zeta$ the complex parameter of this domain. At last, by 
reflecting this domain with respect to the circle $C$ that bounds one of 
these disks, $\Sigma$ can be assumed to be a once punctured bounded domain 
$\Omega \setminus \{z_0\}\subset \C$ whose boundary is made up by
$n$ disjoint circles, the exterior one being precisely $C$.

Once here, by the representation formula we obtain that the parabolic 
affine sphere $\cS_{\phi}=(x,y,\phi(x,y))$ can be conformally parametrized 
as $\psi : \Omega\setminus \{z_0\} \flecha \R^3$, where $\psi$ is given by 
\eqref{rep} in terms of two holomorphic functions $F,G:\Omega\setminus 
\{z_0\} \flecha \C$. Moreover, the map $G +\overline{F} :\Omega\setminus 
\{z_0\} \flecha \R^2 \setminus \{p_1, \dots , p_n\} $ is a global 
diffeomorphism, $G +\overline{F}$ is constant along each boundary circle 
in $\overline{\Omega}$, and $G$ has a pole of order one at $z_0$ 
\cite{FMM2}. 

We recall at this point another result in \cite{FMM2}, which describes the 
asymptotic behaviour of parabolic affine spheres at infinity. This result 
tells in our situation that the solution $\phi(x,y)$ to \eqref{monge} we 
started with has the decomposition 
  \begin{equation}\label{elipse}
\phi (x,y)= E_{\phi} (x,y) + \alfa \log |\zeta|^2 + O(1). 
  \end{equation}
Here $\zeta$ is the conformal parameter of the above mentioned complex 
domain consisting of $\C$ with $n$ disks removed, $\alfa\in \R$, $O(1)$ 
stands for a term bounded in absolute value by a constant, and $E_{\phi} 
(x,y)$ is a quadratic polynomial. When 
$k>0$ is a large positive number, the ellipse $E_{\phi} (x,y)=k$ describes the shape of
$\phi (x,y)$ at infinity, and is called the \emph{ellipse at infinity} of
$\phi(x,y)$.

Now let $\Phi \in \mathcal{E}$ be an equiaffine transformation as in 
\eqref{equi}, and let $\phi'(x',y')$ be the solution to \eqref{monge} 
associated to the embedded parabolic affine sphere $\mathcal{S}' = \Phi 
\circ \cS_{\phi} $. Then $\phi' (x',y'):\R^2 \setminus \{q_1, \dots, q_n\} 
\flecha \R$ is a global solution to \eqref{monge}, where the isolated 
singularities $q_1, \dots, q_n$ are in general different from $p_1,\dots, 
p_n$. Thus, $\phi'(x,',y')$ has all the properties established for 
$\phi(x,y)$, and it becomes clear from the decomposition \eqref{elipse}
that its ellipse at infinity differs from the one of $\phi$ only by an 
affine transformation of $\R^2$. With this, it is easy to choose an 
adequate equiaffine transformation $\Phi$ so that the ellipse at infinity 
of $\phi'(x',y')$ is actually a circle. 

Consider next the holomorphic data of the embedded parabolic affine sphere 
$\mathcal{S}'$, denoted by $G^*, F^*$, and defined on
$\Omega \setminus \{z_0\} $. Then $G^*$ has a simple pole at $z_0$ and, as
the ellipse at infinity of $\mathcal{S}'$ is a circle, $F^*$ extends 
holomorphically about $z_0 $ (see \cite{FMM2}). 

Furthermore, by performing an adequate dilatation and rotation of 
$\Omega$ if necessary, we may assume that the residue of $G^*$ at $z_0$ is
$2$.

With all of this we get that $$\psi_1 = \frac{1}{2} \Re \left( G^*+ F^* 
\right) :\Omega \setminus \{z_0\} \flecha \R$$ is harmonic, constant on 
each boundary component of $\Omega$, and $(G^*+F^*)/2$ has a simple pole 
of residue one at $z_0$. Thus, by uniqueness of the holomorphic map $p(z)$ 
defined in Section 2 we infer that $G^*(z) +F^*(z) =2 p(z) + c_1$, where 
$c_1 \in \C$. Analogously, $G^*(z) -F^*(z) = 2 q(z) + c_2$ for $c_2 \in \C$.
Finally we conclude that $$G^*(z)= p(z) +q(z) + d_1, \hspace{1cm} F^*(z)= 
p(z) -q(z) +d_2,$$ for constants $d_1,d_2\in \C$. This lets us conclude 
that the global solution $\phi'(x',y')$, and thus the original global 
solution $\phi(x,y)$, differs from one of the examples in Theorem 
\ref{ejes} only by an equiaffine transformation. This proves the first 
claim. 

The second claim is essentially direct. We assign to a specific once 
punctured complex domain in $\C$ the example constructed in Theorem 
\ref{ejes}. As changing the domain within its conformal equivalence class 
does not change the surface, we have a well defined mapping from conformal 
equivalence classes of once punctured domains with $n$ boundary curves 
into global solutions to \eqref{monge} in the $n$-times punctured plane, 
modulo equiaffine transformations. This mapping is surjective by the first 
part of the theorem, and injective because solutions to \eqref{monge} 
differing just by an equiaffine transformation have the same underlying 
conformal structure. This ends up the proof. 
\end{proof1}

\begin{remark}
The two classical Jörgens theorems indicate that there is exactly one (up 
to equiaffine transformations) global solution to \eqref{monge} in the 
whole plane $\R^2$, and in the once punctured plane. This does not hold 
anymore when we consider the $n>1$ times punctured plane, as shown by 
Corollary \ref{space}. In any case, given $n>1$ and a conformal 
equivalence class of 
$\C$ with $n$ disjoint disks removed, there exists a unique (up to equiaffine
transformations) global solution to \eqref{monge} in an $n$-times 
punctured plane whose underlying conformal structure is the previously 
given one. So, unicity in the $n=1$ case comes from the fact that all 
domains of the form $\C \setminus D$, where $D$ is a disk in $\C$, are 
conformally equivalent. 
\end{remark}

Theorem \ref{main} can also be stated in the geometric context of 
parabolic affine spheres. Indeed, if a parabolic affine sphere with a 
finite number of singularities is the boundary of a convex body then it 
must be regular at infinity (see \cite{FMM2}), and reasoning as in Theorem 
\ref{ejes} we see that it must be an entire graph. So, we have: 

\begin{cor}\label{Pcv}
All embedded parabolic affine spheres in $\R^3$ that have a finite number 
of isolated singularities and are the boundary of a convex body of 
$\R^3$ must be equiaffinely equivalent to one of the examples in Theorem
\ref{ejes}. 
\end{cor}
Actually, the condition of being the boundary of a convex body can be 
substituted by the property of being closed in $\R^3$. Without going into 
detail, we just indicate that this fact can be proved as follows: given a 
closed parabolic affine sphere $S$, its Legendre transform surface must 
lie on the boundary of a convex body, and so $S$ must be globally convex. 
Thus, Corollary \ref{Pcv} can be applied. 

\section{Explicit construction for one and two punctures}

To begin this last Section, we show that Theorem \ref{main} for $n=1$ 
recovers the classical description by Jörgens of the global solutions to 
\eqref{monge} in the once punctured plane. To do so, we just need to 
analyze the example with $n=1$ in Theorem \ref{ejes}. Obviously, we may 
assume there that $\Omega = \D$, the unit circle, and $z_0 =0$. 

Let $z_0 \in \D$. The Green function in $\D$ with respect to $z_0$ is then 
$$\cG (z,z_0)= \log \left|\frac{1- z \overline{z_0}}{z -
z_0}\right| .$$ Deriving this expression with respect to $z_0$ as in 
\eqref{aru}, \eqref{ardo} we get the harmonic functions $$ u(z,z_0)=\Re 
\left( \frac{1}{z-z_0} - \frac{z}{1 - z\overline{z_0}}\right), 
\hspace{0.6cm} v(z,z_0)=\Im \left( \frac{1}{z-z_0} + \frac{z}{1 - 
z\overline{z_0}}\right) .$$ Therefore, letting $z_0=0$ we see that the 
canonical functions $p(z),q(z)$ in 
$\D^*$ are $$p(z)= \frac{1}{z} -z, \hspace{1cm}
q(z)= \frac{1}{z} +z.$$ Finally, by \eqref{fg} the holomorphic data of the 
solution to \eqref{monge} are $G(z)=2 /z$ and $F(z)=-2z$. These are known 
to be the holomorphic data (up to the conformal change $w=2/z$) of the 
rotational example in the introduction (see \cite{FMM1}). Also observe 
that $(G + \overline{F})|_{\S^1} =0$. Hence, Jörgens' theorem is 
recovered. 

\begin{figure}[h]\label{rev}
  \begin{center}
\includegraphics[clip,width=5.5cm]{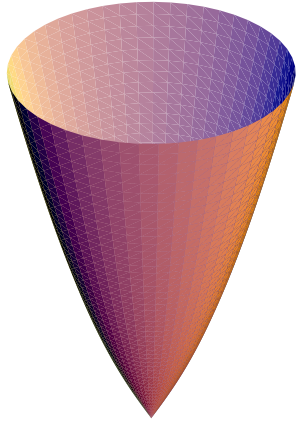}\vspace{-2cm}
\caption{The unique solution to the Monge-Ampère equation in the once 
punctured plane, modulo equiaffine transformations.} 
 \end{center}
\end{figure}

It is also possible to write down explicitly in terms of theta functions 
all global solutions to \eqref{monge} in a twice punctured plane 
$\R^2 \setminus \{p_1,p_2\} $, as we show next.

 \begin{figure}[h]
  \begin{center}
    \begin{tabular}{cc}
\includegraphics[clip,width=7cm]{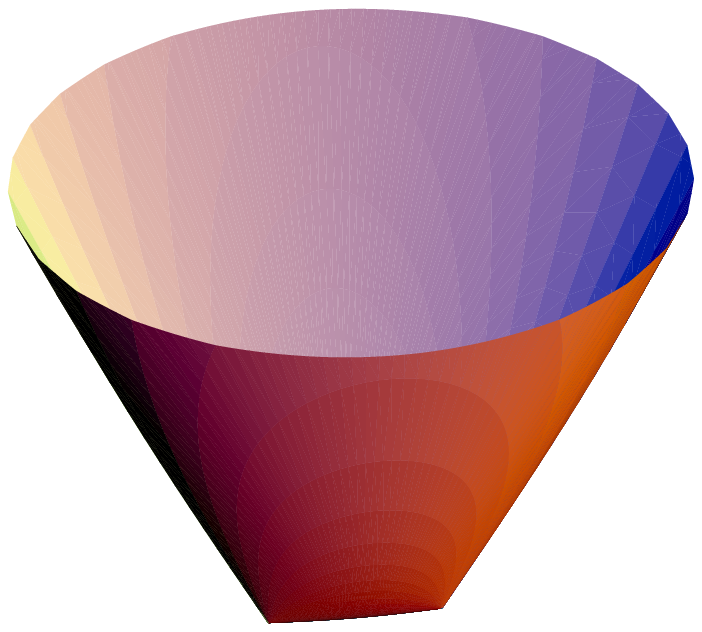} &
\includegraphics[clip,width=6.5cm]{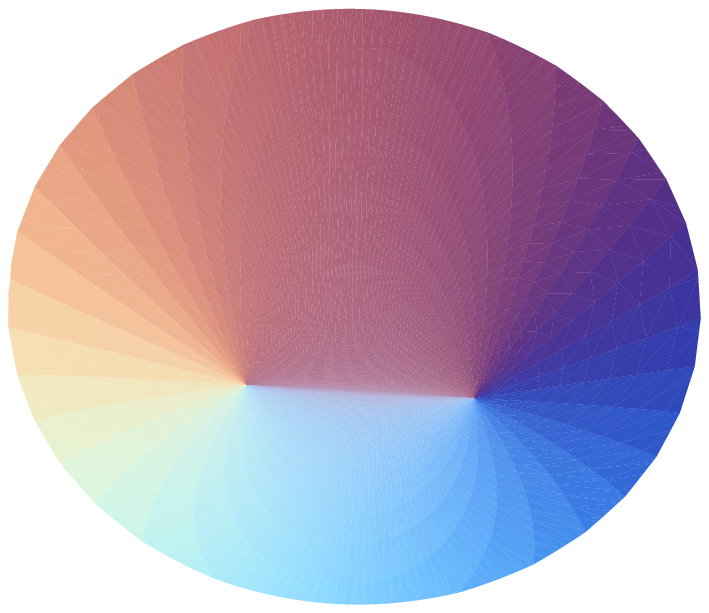}\vspace{-.5cm}
\end{tabular}
\caption{Two different views of a solution to the Monge-Ampère equation in 
a twice punctured plane.} 
 \end{center}
\end{figure}

Let $\theta_0 (w)$ be the Jacobi theta function (see \cite{CoHi} for 
instance) 
$$\theta_0(w) =C \prod_{k=1}^{\8} \left(1 - r^{2k-1} e^{2 \pi i w}\right)
\left(1 - r^{2k-1} e^{-2 \pi i w}\right),$$ where $C=\prod_{k=1}^{\8} 
(1-r^{2k})$ and $0<r<1$. Consider also $\mathbb{A}_r$ to be the annulus 
given by $\mathbb{A}_r = \{z : r<|z| <1\}$, and let $\vart_1 (z)$ denote 
the \emph{annular} Jacobi theta function given by $\vart_1 (e^{2\pi i w})= 
\theta_0 (r w)$. In other words, 
 \begin{equation}\label{tet}
 \vart_1 (z)= C \left( 1- \frac{1}{z}\right)
\prod_{k=1}^{\8} \left( 1- r^{2k} z \right) \left( 1- r^{2k}/ z \right). 
 \end{equation}
It comes then clear from this expression that 
 \begin{equation}\label{tetrel}
 \vart_1 (z)= \overline{\vart_1(\bar{z})} = -r^2 z \, \vart_1 (r^2 z),
 \hspace{0.7cm} \vart_1 (z /r^2)= -z \vart_1 (z),
 \end{equation}
and thus that 
 \begin{equation}\label{tetre2}
 \vart_1 '(z) = - r^2 \vart_1 (r^2 z) - r^4 z \vart_1'(r^2 z),
 \hspace{0.5cm} \vart_1'(z/r^2)= -r^2 \vart_1(z) - r^2 z \vart_1 '(z).
 \end{equation}

A formula for the Green function of an annulus of the type 
$\{r^{1/2} < |z|< r^{-1/2}\}$ for $r \in (0,1)$ can be found in
\cite[pg. 386-387]{CoHi}. By the conformal invariance property of the 
Green function, we obtain from this formula that the Green function 
$\cG(z,z_0)$ of $\mathbb{A}_r$ with respect to a point $z_0\in \mathbb{A}_r$ can be
obtained in terms of 
$\vart_1$ as $$\cG (z,z_0)= \log |z| \left( 1 +\frac{\log
|z_0|}{\log r}\right) -\log \left|\frac{\vart_1 ( z_0 /z)}{\vart_1 (z 
\overline{z_0})} \right| . 
$$ In other words, $$\cG(z,z_0)= \Re h(z) = \Re \left\{\log
\left( \frac{z^{1 +\log |z_0| / \log r}}{\vart_1 (z_0/z) / \vart_1 
(\overline{z_0} z)}\right) \right\} .$$ At this point we can find the 
harmonic function in \eqref{aru} by differentiation. We obtain $$u (z)= 
\Re \left\{ \frac{\log z}{\log r} \left(\frac{1}{z_0} + 
\frac{1}{\overline{z_0}}\right)\right\} -\Re \left \{ \frac{\vart_1'(z_0 
/z)}{z \vart_1 (z_0 /z)} \frac{\vart_1'(\overline{z_0} z) z}{\vart_1 
(\overline{z_0}z)} \right\} .$$ Noting that $$ \frac{\vart_1'(z_0 /z)}{z 
\vart_1 (z_0 /z)} -\frac{\vart_1'(\overline{z_0} z) z}{\vart_1 
(\overline{z_0}z)}$$ is a well defined meromorphic function on 
$\mathbb{A}_r$, the way we defined the holomorphic function $p(z)$ in
\eqref{truqui} suggests that 
 \begin{equation}\label{p2}
 p(z)=\frac{\vart_1'(\overline{z_0}
z) z}{\vart_1 (\overline{z_0}z)} - \frac{\vart_1'(z_0 /z)}{z \vart_1 (z_0 
/z)} . 
 \end{equation}
The validity of this assertion is easily checked: $p(z)$ is holomorphic on 
$\mathbb{A}_r \setminus \{z_0\} $, it has a simple pole of residue $1$ at
$z_0$, and using \eqref{tetrel} and \eqref{tetre2} we see that
 \begin{equation}\label{cambi}
p (1/\bar{z}) + \overline{p(z)} =0, \hspace{1cm} p (r^2 /\bar{z}) + 
\overline{p(z)} = -\frac{1}{z_0} - \frac{1}{\overline{z_0}}. 
 \end{equation}
So, the real part of $p$ is constant on the boundary curves $|z|=r,1$. 

The corresponding process starting this time with the harmonic function 
defined via \eqref{ardo} lets us recover the other canonical holomorphic 
function $q(z):\mathbb{A}_r \setminus \{z_0\} \flecha \C$ as 
  \begin{equation}\label{p3}
 q(z)=-\frac{\vart_1'(z_0 /z)}{z \vart_1 (z_0 /z)} -\frac{\vart_1'(\overline{z_0}
z) z}{\vart_1 (\overline{z_0}z)}. 
 \end{equation}
Putting together \eqref{p2} and \eqref{p3}, and recalling the formula 
\eqref{fg}, we obtain the holomorphic data on $\mathbb{A}_r \setminus 
\{z_0\} $ 

 \begin{equation}\label{data2}
 G(z)= -\frac{2 \vart_1'(z_0 /z)}{z \vart_1 (z_0 /z)}, \hspace{0.6cm} F(z)= \frac{2\vart_1'(\overline{z_0}
z) z}{\vart_1 (\overline{z_0}z)}. 
 \end{equation}

Finally, by Theorem \ref{main}, all solutions to \eqref{monge} in a twice 
punctured plane are recovered, up to equiaffine transformations, by the 
holomorphic data in \eqref{data2}. The 2-parametric family of such 
solutions described by Corollary \ref{space} is given by the variation of 
$r\in (0,1)$ and $|z_0| \in (r,1)$. In the solution to \eqref{monge} we have just
constructed, the singularities are located at the points $(0,0)$ and $(\Re 
(-2 /\overline{z_0}), \Im (-2/\overline{z_0}))$ of $\R^2$. This follows 
from \eqref{rep} and the relation \eqref{cambi}.

\def\refname{References}

\end{document}

\end{document}